\documentclass[11pt]{article}
\usepackage[a4paper,margin=25mm]{geometry}
\usepackage[T1]{fontenc}
\usepackage{lmodern,microtype,amsmath,amssymb,amsthm,mathtools}
\usepackage[hidelinks]{hyperref}
\newtheorem{theorem}{Theorem}
\newtheorem{lemma}[theorem]{Lemma}

\newtheorem{proposition}[theorem]{Proposition}
\newcommand{\dd}{\,d}
\newcommand{\ii}{\mathrm i}

\DeclareMathOperator{\Impart}{Im}
\DeclareMathOperator{\li}{li}
\author{Martin Belton}
\date{6 November 2026}
\usepackage{xcolor,graphicx}
\hypersetup{pdftitle={Riemann functions for almost primes and their counting errors},pdfauthor={Martin Belton}}
\title{Riemann functions for almost primes and their counting errors}
\begin{document}
\maketitle
\begin{abstract}
The real-axis jump of the prime-zeta function gives an integral form of
Gram's series for Riemann's $R$-function. Applying the same operation to
the almost-prime Dirichlet series defines functions $R_k$, with $R_1=R$.
We prove absolute convergence, recover the classical fixed-order
Selberg--Delange expansion, and identify the complementary Perron
contribution to the exact count. A convergent incomplete-gamma endpoint
expansion, together with the remaining lower cut, gives an exact
representation of $R_k$, with explicit formulas for ranks one through
four and a coefficient recurrence for every fixed rank. We then study
the errors $R_k-N_k$.
Under RH their normalized mean-square zero spectra retain multiplicities;
with simple zeros the leading oscillations agree for every fixed $k$.
An additional weighted logarithmic-gap moment gives a second-order
spectrum and a positive constant in the $(\log\log X)^{-2}$ mean-square
cancellation rate for consecutive values of $k$. The global spectral
passages are justified by a uniform finite contour and a
Gaussian--Laplace block estimate.
\end{abstract}

\section{Introduction}
Let $\Omega(n)$ count prime factors with multiplicity, and write
\[
 N_k(x)=\#\{n\le x:\Omega(n)=k\}.
\]
Thus $N_1$ counts primes, $N_2$ counts semiprimes, and $N_3$ counts
products of three primes, including repeated factors. We ask how
Riemann's smooth prime-counting function $R$ extends to these counts,
and what can be said about the resulting arithmetic errors.

The starting point is a prescribed integral over the complete positive
real cut of the almost-prime Dirichlet series. Its normalization gives
$R_1=R$ exactly. At every fixed inverse-logarithmic order its expansion
is the classical Selberg--Delange expansion. The definition also retains
the smaller real-axis contributions that such an expansion does not
determine. A finite Perron identity separates this smooth contribution
from the complementary contour terms.

The construction, convergence and fixed-order expansion are
unconditional. The error laws have successive explicit hypotheses:
RH gives the multiplicity-sensitive mean-square zero expansion;
simplicity aligns the leading spectra for consecutive ranks; an
additional logarithmic-gap moment gives the second-order cancellation
rate. All asymptotic statements keep $k$ fixed. The laws are
logarithmically averaged, not pointwise approximations or estimates
proving RH. The uniform contour and Gaussian--Laplace lemmas below
supply the passage from local zero contributions to the actual errors.

The ingredients of the construction are classical
\cite{Froberg,Mathar,SD}; explicit semiprime asymptotics are given by
Cri\c{s}an and Erban \cite{CrisanErban}. No new Selberg--Delange
coefficients are claimed. Meng \cite{Meng} proves related alternating
factorial/logarithmic scaling for differences of almost-prime counts
in residue classes. Here the principal zeta function and its complete
real cut are retained, with a multiplicity defect and a second-order
error spectrum. This paper contains the construction and spectral
proofs; coefficient predictors and their relation to Liouville parity
are developed separately in the companion manuscript \cite{BeltonParity}.

A related real-axis contour construction appears in Lucia's
MathOverflow answer \cite{LuciaMO}, which credits the underlying method
to Selberg. It obtains a slit-integral main term for the coefficients of
$(\log\zeta(s))^k/k!$ and distinguishes those weighted coefficients from
exact almost-prime counts. Here we use the full series $F_k$, prescribe
its complete positive real cut, and normalize by $R_1=R$; the
real-axis slit-contour principle itself is not claimed as new.

\section{From Gram's series to almost primes}\label{sec:definition}
Put $B(t)=\sum_{j\le t}\mu(j)/j$, with $B(t)=0$ for $t<1$.
Gram's series has the integral form
\begin{equation}\label{eq:gram}
 R(x)=1+\sum_{n\ge1}\frac{(\log x)^n}{n\,n!\,\zeta(n+1)}
     =\int_0^1 B(1/\sigma)\frac{x^\sigma}{\sigma}\dd\sigma,
 \qquad x>0.
\end{equation}
Indeed the quantitative prime number theorem gives
$B(t)\ll_D(\log t)^{-D}$ for every fixed $D>0$ \cite{MV},
so the integral converges absolutely. Partial summation gives
\[
 \int_0^1B(1/\sigma)\sigma^{v-1}\dd\sigma
       =\frac1{v\zeta(1+v)}\quad(v>0).
\]
Its limit at $v=0$ is $1$. Expanding $x^\sigma$ under the integral
therefore recovers each coefficient in Gram's series.

For $\Re s>1$, let $P(s)=\sum_p p^{-s}$. Then
\begin{equation}\label{eq:generator}
 P(s)=\sum_{j\ge1}\frac{\mu(j)}j\log\zeta(js),\qquad
 F_k(s):=\sum_{\Omega(n)=k}n^{-s}
 =[z^k]\exp\left(\sum_{j\ge1}\frac{z^jP(js)}j\right).
\end{equation}
Only $j\le k$ contributes to the coefficient. Define the polynomial
$H_k$ by $F_k(s)=H_k(P(s),\ldots,P(ks))$. Explicitly,
\[
 H_k(z_1,\ldots,z_k)=
 \sum_{\nu_1+2\nu_2+\cdots+k\nu_k=k}
          \prod_{j=1}^k\frac{z_j^{\nu_j}}{j^{\nu_j}\nu_j!}.
\]
Away from reciprocal integers, prescribe the lower bank of the real cut by
\[
 P_-(\sigma)=\sum_{j\ge1}\frac{\mu(j)}j\log|\zeta(j\sigma)|
                   +\ii\pi B(1/\sigma),\qquad \sigma>0.
\]
This is continuation below the real poles, where
$\log\zeta(t-\ii0)=\log|\zeta(t)|+\ii\pi$ for $0<t<1$.
The other bank is its conjugate. The series converges locally normally
away from reciprocal endpoints because its tail decays geometrically.
Replacing the prime series by the almost-prime series in the same
jump operation defines
\begin{equation}\label{eq:Rk}
 \Delta_k(\sigma)=\pi^{-1}\Impart
 H_k(P_-(\sigma),\ldots,P_-(k\sigma)),\qquad
 R_k(x)=\int_0^1\Delta_k(\sigma)\frac{x^\sigma}{\sigma}\dd\sigma.
\end{equation}

The linear term $P(k\sigma)/k$ in $H_k$ contributes
$R(x^{1/k})/k$. For example, writing
$U(\sigma)=\sum_{j\ge1}\mu(j)\log|\zeta(j\sigma)|/j$, we obtain
\begin{equation}\label{eq:semiprime}
 R_2(x)=\tfrac12R(\sqrt{x})
       +\int_0^1 B(1/\sigma)U(\sigma)\frac{x^\sigma}{\sigma}\dd\sigma.
\end{equation}
Here $H_2(w_1,w_2)=(w_1^2+w_2)/2$; the $w_2$ term corrects
the prime-square diagonal. Similarly
$H_3=(w_1^3+3w_1w_2+2w_3)/6$ includes $p^3$ and $p^2q$.
The reciprocal endpoints are improper endpoints, not numerical cutoffs.

\begin{lemma}\label{lem:cut}
For fixed $k$, $\int_0^1|\Delta_k(\sigma)|\dd\sigma/\sigma<\infty$.
Thus $R_k(e^w)$ is entire in $w$, and the portion $0<\sigma<a$ is
$O_{k,a}(x^a)$ for $x\ge1$. In particular, \eqref{eq:Rk} gives $R_1=R$ by \eqref{eq:gram}.
\end{lemma}
\begin{proof}
The quantitative prime number theorem and partial summation give
$B(t)\ll_D(\log t)^{-D}$ for every fixed $D>0$ \cite{MV}.
Let $I_q=(1/(q+1),1/q)$ and $\dd\nu_q=q\dd\sigma/\sigma$.
For fixed $m,p\le k$,
\begin{equation}\label{eq:cell}
 \|P_-(m\sigma)\|_{L^p(\nu_q)}\ll_k\log(q+2),\qquad
 |B(1/(m\sigma))|\ll_{k,D}(\log q)^{-D}.
\end{equation}
Here is a direct bound for the first assertion. Set $v=1/(m\sigma)$
and split the defining real sum at $v/2$ and $2v$. For $j<v/2$,
$\log|\zeta(j/v)|$ is bounded, giving $O(\log v)$ by absolute
summation. For $j>2v$, exponential decay gives $O(1)$. In the middle
range use $\log|\zeta(t)|=-\log|t-1|+O(1)$.
After removing the two indices nearest $v$, the sum is $O(1)$,
by comparison with $v^{-1}\sum_{h\le2v}\log(2v/h)$.
The two remaining terms have $L^p$ norm $O_k((1+\log v)/v)$: the interval
$q/m<v<(q+1)/m$ contains no integer, and
$\int_0^1|\log t|^p\dd t=p!$. This proves~\eqref{eq:cell}.

For a product of $r\le k$ factors, telescope its imaginary part
against the product of their real parts. Each term then contains one
factor $B(1/(m\sigma))$. H\"older and~\eqref{eq:cell} show that its
integral on $I_q$ is
$O_{k,D}(q^{-1}(\log q)^{r-1-D})$. Choose $D=k+2$ and sum over $q$.
The finitely many remaining cells have only integrable endpoint
logarithms. Dominated differentiation proves the entire extension and
the lower-cut bound.

\end{proof}

\section{The classical expansion and the full real cut}\label{sec:selberg}
Initially for $\Re s>1$ and $|z|$ small, the generating function factors as
\begin{equation}\label{eq:factor}
 \mathcal F(s,z)=\prod_p(1-zp^{-s})^{-1}
       =\zeta(s)^zG(s,z),\qquad
 G(s,z)=\exp\left(\sum_{j\ge2}\frac{z^j-z}{j}P(js)\right).
\end{equation}
The function $G$ is analytic near $(s,z)=(1,0)$.
Define the local coefficient generator
\begin{equation}\label{eq:bj}
 C(t,z)=\frac{(t\zeta(1+t))^zG(1+t,z)}{1+t}
        =\sum_{j\ge0}b_j(z)t^j.
\end{equation}
The logarithm of $t\zeta(1+t)$ is normalized to vanish at $t=0$.
Equivalently, the same coefficients are defined by
\begin{equation}\label{eq:bj-primezeta}
 C(t,z)=\frac{1}{1+t}
 \exp\left(z\{P(1+t)+\log t\}
       +\sum_{m\ge2}\frac{z^mP(m(1+t))}{m}\right).
\end{equation}
Here $P(1+t)+\log t$ denotes its analytic continuation at $t=0$,
not two separately evaluated singular terms: it equals
$\log(t\zeta(1+t))-\sum_{m\ge2}P(m(1+t))/m$.
Identity \eqref{eq:bj-primezeta} follows by substituting
$\log\zeta(1+t)=P(1+t)+\sum_{m\ge2}P(m(1+t))/m$
into \eqref{eq:bj}.

\begin{proposition}\label{prop:SD}
For every fixed $k\ge1$ and $J\ge0$, as $x\to\infty$,
\begin{align}\label{eq:SD}
 R_k(x)&=x\sum_{j=0}^J\frac1{(\log x)^{j+1}}
            [z^k]\frac{(\log x)^z b_j(z)}{\Gamma(z-j)}\\
 &\quad+O_{k,J}\left(
       \frac{x(1+\log\log x)^{k-1}}{(\log x)^{J+2}}\right).\notag
\end{align}
These are the classical Selberg--Delange coefficients.
In particular
\[
 R_k(x)\sim\frac{x(\log\log x)^{k-1}}{(k-1)!\log x}\sim N_k(x).
\]
\end{proposition}
\begin{proof}
Choose $0<\delta<1/2$ and $0<r<1$ so that $C$ is analytic on
a neighbourhood of $|t|\le\delta$, $|z|\le r$.
On the upper part of the real cut, put $\sigma=1-v$.
The boundary difference in \eqref{eq:factor} gives
\[
 \sum_{k\ge1}\frac{\Delta_k(1-v)}{1-v}z^k
        =\frac{\sin\pi z}{\pi}v^{-z}C(-v,z),
 \qquad 0<v<\delta.
\]
For a Taylor term, Euler's reflection identity yields
\[
 \frac{(-1)^j\sin\pi z}{\pi}
 \int_0^\infty e^{-v\log x}v^{j-z}\dd v
       =\frac{(\log x)^{z-j-1}}{\Gamma(z-j)}.
\]
Write $C(t,z)=\sum_{j=0}^Jb_j(z)t^j+t^{J+1}E_J(t,z)$.
The finitely many needed $z$ derivatives of $E_J$ are uniformly
bounded on a smaller polydisc, so
\[
 [z^k]\left\{\frac{\sin\pi z}{\pi}v^{-z}(-v)^{J+1}E_J(-v,z)\right\}
       =O_{k,J}\bigl(v^{J+1}(1+|\log v|)^{k-1}\bigr).
\]
The degree is $k-1$ because $\sin\pi z$ vanishes at zero.
Laplace scaling proves the error in \eqref{eq:SD}. Extending each
retained integral to infinity, and omitting the rest of the real
cut, cost $O_{k,J,\delta}(x^{1-\delta})$ by Lemma~\ref{lem:cut}.
The factorization and coefficient rule are precisely the local
Selberg--Delange construction \cite{SD,MV}.
Finally $b_0(z)=G(1,z)$, $G(1,0)=1$ and
$1/\Gamma(z)=z+O(z^2)$ give the leading term.
\end{proof}

\paragraph{An exact convergent endpoint expansion.}
The same coefficients also give a convergent representation when the
endpoint integral is kept on a finite interval rather than extended
to infinity. This is distinct from summing the asymptotic series
\eqref{eq:SD} at a fixed value of $x$.

\begin{proposition}[Convergent endpoint expansion]\label{prop:endpoint}
Choose $0<\delta<R<1/2$ and $0<r<1$ such that $C(t,z)$ is analytic
on a neighbourhood of $|t|\le R$, $|z|\le r$. For $k\ge1$ and $x>1$,
put
\begin{equation}\label{eq:lowercut-delta}
 T_{k,\delta}(x)=\int_0^{1-\delta}
             \Delta_k(\sigma)\frac{x^\sigma}{\sigma}\dd\sigma.
\end{equation}
Then
\begin{equation}\label{eq:exact-endpoint}
{\displaystyle
 R_k(x)=T_{k,\delta}(x)
  +x[z^k]\frac{\sin(\pi z)}{\pi}
    \sum_{j=0}^{\infty}(-1)^j b_j(z)
      \frac{\gamma(j+1-z,\delta\log x)}{(\log x)^{j+1-z}} .}
\end{equation}
Here $\gamma(a,v)=\int_0^v e^{-u}u^{a-1}\dd u$ is the lower
incomplete gamma function, with $\Re a>0$ and $v>0$ \cite[\S8.2]{DLMF}.
The series before coefficient extraction converges absolutely and
uniformly for $|z|\le r$, locally uniformly for $x>1$.
Coefficient extraction may therefore be interchanged with the sum.

More explicitly, let $U_{k,\delta,J}(x)$ denote the upper-cut sum in
\eqref{eq:exact-endpoint} truncated at $j=J$, and set
\[
 M_{R,r}=\max_{\substack{|t|=R\\|z|=r}}|C(t,z)|,
 \qquad S_r=\max_{|z|=r}\frac{|\sin(\pi z)|}{\pi}.
\]
For every $J\ge0$,
\begin{align}\label{eq:endpoint-geometric-bound}
 &|R_k(x)-T_{k,\delta}(x)-U_{k,\delta,J}(x)|\notag\\
 &\qquad\le
 \frac{x\,r^{-k}S_rM_{R,r}\,\delta^{1-r}}{1-\delta/R}
 \frac{(\delta/R)^{J+1}}{J+2-r}.
\end{align}
\end{proposition}
\begin{proof}
The upper-cut boundary identity used in the proof of
Proposition~\ref{prop:SD} gives
\[
 R_k(x)-T_{k,\delta}(x)
 =x[z^k]\frac{\sin(\pi z)}{\pi}
       \int_0^\delta e^{-v\log x}v^{-z}C(-v,z)\dd v.
\]
Cauchy's estimate gives $|b_j(z)|\le M_{R,r}R^{-j}$ on $|z|\le r$.
Since $\delta<R$ and $r<1$, the integrated Taylor series is dominated
by a constant times $v^{-r}/(1-v/R)$, integrable on $(0,\delta)$.
Termwise integration and coefficient extraction are justified by this
bound, using a $z$-circle of radius $r$. Each integral is exactly
\[
 \int_0^\delta e^{-v\log x}v^{j-z}\dd v
  =\frac{\gamma(j+1-z,\delta\log x)}{(\log x)^{j+1-z}},
\]
which proves \eqref{eq:exact-endpoint}. The same domination is uniform
on compact subsets of $x>1$.

For the truncation bound, the Taylor remainder on $0\le v\le\delta$
is at most
\[
 \frac{M_{R,r}}{1-\delta/R}(v/R)^{J+1}.
\]
Apply Cauchy's coefficient bound $|[z^k]f|\le r^{-k}\max_{|z|=r}|f|$,
use $|v^{-z}|\le v^{-r}$ and $e^{-v\log x}\le1$, and integrate
$v^{J+1-r}$ from $0$ to $\delta$. This gives
\eqref{eq:endpoint-geometric-bound}.
\end{proof}

Replacing the finite integral by an integral to infinity changes
$\gamma(j+1-z,\delta\log x)$ to $\Gamma(j+1-z)$; reflection then
gives the coefficients $1/\Gamma(z-j)$ in \eqref{eq:SD}.
That passage is valid at each fixed truncation order in the asymptotic
argument, not as an unrestricted exchange with the infinite $j$-sum.
The incomplete-gamma series is the exact upper-cut contribution;
only after adding $T_{k,\delta}$ is it the full prescribed $R_k$.
The convergence bound does not by itself certify a numerical
implementation: the coefficient bounds and the lower-cut integral
must also be evaluated or enclosed. Nor does it estimate $N_k-R_k$.

\subsection{Regular prime-zeta coefficients}
The following formulas make Proposition~\ref{prop:endpoint} explicit.
Retain its choices of $\delta,R,r$, and put
\begin{equation}\label{eq:endpoint-h}
 h(t)=P(1+t)+\log t=\sum_{j\ge0}h_jt^j.
\end{equation}
The cancellation at $t=0$ is understood analytically. Equivalently,
\begin{equation}\label{eq:endpoint-h-regular}
 {\displaystyle
 h(t)=\log\!\bigl(t\zeta(1+t)\bigr)
      -\sum_{m=2}^{\infty}\frac{P(m(1+t))}{m}.}
\end{equation}
Normalize the first logarithm to vanish at $t=0$. Since $R<1/2$,
the prime-zeta sum is normally convergent on $|t|\le R$; its
termwise derivatives of every fixed order also converge there.
It follows that
\begin{equation}\label{eq:endpoint-hj}
 {\displaystyle
 h_j=[t^j]\log\!\bigl(t\zeta(1+t)\bigr)
       -\frac1{j!}\sum_{m=2}^{\infty}m^{j-1}P^{(j)}(m).}
\end{equation}
This includes $j=0$, with $P^{(0)}=P$. Here and below $[t^j]$
denotes the Taylor coefficient, not the $j$th derivative.
In particular the factorial division in \eqref{eq:endpoint-hj}
is already included.

Let $M=\gamma_E-\sum_{m\ge2}P(m)/m$ be the prime
Meissel--Mertens constant. We use the Stieltjes convention
\cite[\S25.2]{DLMF}
\[
 \zeta(1+t)=\frac1t+\sum_{n\ge0}\frac{(-1)^n\gamma_n}{n!}t^n,
 \qquad \gamma_0=\gamma_E.
\]
Expansion of the logarithm gives
\begin{align}\label{eq:endpoint-h-first}
 h_0={}&M-\gamma_E,\notag\\
 h_1={}&\gamma_E-\sum_{m\ge2}P'(m),\notag\\
 h_2={}&-\gamma_1-\frac{\gamma_E^2}{2}
            -\frac12\sum_{m\ge2}mP''(m),\notag\\
 h_3={}&\frac{\gamma_2}{2}+\gamma_E\gamma_1+
           \frac{\gamma_E^3}{3}
            -\frac16\sum_{m\ge2}m^2P'''(m),\notag\\
 h_4={}&-\frac{\gamma_3}{6}-\frac{\gamma_E\gamma_2}{2}
          -\frac{\gamma_1^2}{2}-\gamma_E^2\gamma_1-
           \frac{\gamma_E^4}{4}\notag\\
       &\hspace{8mm}-\frac1{24}\sum_{m\ge2}m^3P^{(4)}(m).
\end{align}
These first coefficients illustrate the rule. The exact infinite
endpoint expansions below use all required $h_j$, not only
$h_0,\ldots,h_4$.

\subsection{Incomplete-gamma moments and the first four ranks}
For integers $j,\ell\ge0$, $x>1$ and the same $\delta$, define
\begin{equation}\label{eq:ig-moments}
 {\displaystyle
 I_{j,\ell}(x,\delta)=
 \left.\frac{\partial^{\ell}}{\partial z^{\ell}}
 \left[\frac{\gamma(j+1-z,\delta\log x)}
            {(\log x)^{j+1-z}}\right]\right|_{z=0}.}
\end{equation}
Equivalently,
\begin{equation}\label{eq:ig-moments-integral}
 I_{j,\ell}(x,\delta)=\int_0^\delta
       e^{-v\log x}v^j(-\log v)^{\ell}\dd v.
\end{equation}
Differentiation under the integral is justified by
$v^{j-r}(1+|\log v|)^{\ell}$ for $|z|\le r<1$.
Thus the $I_{j,\ell}$ are specified special-function moments,
not additional undetermined constants.

Define the endpoint coefficient sequences
\begin{align}\label{eq:endpoint-ABC}
 A_j&=[t^j]\frac{h(t)}{1+t},\notag\\
 B_j&=[t^j]\frac{h(t)^2+P(2(1+t))}{2(1+t)},\\
 C_j&=[t^j]\frac{h(t)^3/6+h(t)P(2(1+t))/2+P(3(1+t))/3}{1+t}.
 \notag
\end{align}
The notation $A_j$ here denotes a coefficient, as opposed to the
normalized counting-error function $A_k(y)$ introduced later.
The sequences can be evaluated from the $h_j$ and ordinary
prime-zeta derivatives:
\begin{align}\label{eq:endpoint-ABC-convolutions}
 A_j={}&\sum_{n=0}^j(-1)^{j-n}h_n,\notag\\
 B_j={}&\frac12\sum_{n=0}^j(-1)^{j-n}
       \left(\sum_{a=0}^n h_a h_{n-a}
                   +\frac{2^n}{n!}P^{(n)}(2)\right),\notag\\
 C_j={}&\sum_{n=0}^j(-1)^{j-n}\Biggl(
       \frac16\sum_{\substack{a,b,c\ge0\\a+b+c=n}}h_a h_b h_c\notag\\
     &\hspace{19mm}+\frac12\sum_{a=0}^n h_a
                  \frac{2^{n-a}}{(n-a)!}P^{(n-a)}(2)
                     +\frac{3^n}{3n!}P^{(n)}(3)\Biggr).
\end{align}

In the next four formulas write $I_{j,\ell}=I_{j,\ell}(x,\delta)$
to avoid repeating the arguments. Then
\begin{equation}\label{eq:endpoint-R1}
 {\displaystyle
 R_1(x)=T_{1,\delta}(x)+x\sum_{j=0}^{\infty}I_{j,0}.}
\end{equation}
\begin{equation}\label{eq:endpoint-R2}
 {\displaystyle
 R_2(x)=T_{2,\delta}(x)+x\sum_{j=0}^{\infty}
           \left(I_{j,1}+(-1)^jA_jI_{j,0}\right).}
\end{equation}
\begin{equation}\label{eq:endpoint-R3}
 {\begin{aligned}
 R_3(x)=T_{3,\delta}(x)+x\sum_{j=0}^{\infty}\Biggl[
       &\frac{I_{j,2}}2+(-1)^jA_jI_{j,1}\\
       &+\left((-1)^jB_j-\frac{\pi^2}{6}\right)I_{j,0}
                                              \Biggr].
 \end{aligned}}
\end{equation}
\begin{equation}\label{eq:endpoint-R4}
 {\begin{aligned}
 R_4(x)=T_{4,\delta}(x)+x\sum_{j=0}^{\infty}\Biggl[
       &\frac{I_{j,3}}6+\frac{(-1)^jA_j}{2}I_{j,2}\\
       &+\left((-1)^jB_j-\frac{\pi^2}{6}\right)I_{j,1}\\
       &+(-1)^j\left(C_j-\frac{\pi^2}{6}A_j\right)I_{j,0}
                                              \Biggr].
 \end{aligned}}
\end{equation}
All four series converge absolutely for each $x>1$, locally
uniformly in $x$. Each retains the lower-cut integral
$T_{k,\delta}$ from \eqref{eq:lowercut-delta}.
As a check, summing the geometric series in \eqref{eq:ig-moments-integral}
gives
\begin{align}\label{eq:endpoint-R1-check}
 x\sum_{j\ge0}I_{j,0}
  &=x\int_0^\delta\frac{e^{-v\log x}}{1-v}\dd v\notag\\
  &=\li(x)-\li(x^{1-\delta}).
\end{align}
Adding $T_{1,\delta}$, rather than dropping it, completes this
upper-cut contribution to the prescribed Riemann function $R_1=R$.

\subsection{A recurrence for every fixed rank}
Define analytic functions $E_n(t)$ by
\begin{equation}\label{eq:endpoint-E-recurrence}
 E_0(t)=1,\qquad
 E_n(t)=\frac1n\left(h(t)E_{n-1}(t)
       +\sum_{m=2}^{n}P(m(1+t))E_{n-m}(t)\right)
       \quad(n\ge1).
\end{equation}
They are the coefficients of
\begin{equation}\label{eq:endpoint-E-generator}
 \exp\left(zh(t)+\sum_{m\ge2}\frac{z^mP(m(1+t))}{m}\right)
       =\sum_{n\ge0}E_n(t)z^n.
\end{equation}
Differentiation with respect to $z$ proves the recurrence. Put
\begin{equation}\label{eq:endpoint-cjn}
 c_{j,n}=[t^j]\frac{E_n(t)}{1+t}.
\end{equation}
Then $b_j(z)=\sum_{n\ge0}c_{j,n}z^n$, with
\[
 c_{j,0}=(-1)^j,\qquad c_{j,1}=A_j,\qquad
 c_{j,2}=B_j,\qquad c_{j,3}=C_j.
\]
For every fixed $k\ge1$, the general rule is
\begin{equation}\label{eq:endpoint-general-rank}
 {\begin{aligned}
 R_k(x)={}&T_{k,\delta}(x)+x\sum_{j=0}^{\infty}(-1)^j
       \sum_{q=0}^{\lfloor(k-1)/2\rfloor}
                 \frac{(-1)^q\pi^{2q}}{(2q+1)!}\\
       &\hspace{13mm}\times
       \sum_{\ell=0}^{k-2q-1}
          \frac{c_{j,k-2q-1-\ell}}{\ell!}I_{j,\ell}(x,\delta).
 \end{aligned}}
\end{equation}
Only $E_0,\ldots,E_{k-1}$ are needed for rank $k$.
\begin{proof}[Derivation of the rank formulas]
In \eqref{eq:exact-endpoint}, use
\[
 \frac{\sin(\pi z)}{\pi}
       =\sum_{q\ge0}\frac{(-1)^q\pi^{2q}}{(2q+1)!}z^{2q+1},
 \qquad
 \frac{\gamma(j+1-z,\delta\log x)}{(\log x)^{j+1-z}}
       =\sum_{\ell\ge0}\frac{I_{j,\ell}(x,\delta)}{\ell!}z^\ell.
\]
Together with $b_j(z)=\sum_n c_{j,n}z^n$, the condition
$n+\ell+2q+1=k$ gives \eqref{eq:endpoint-general-rank}.
The convergence proved in Proposition~\ref{prop:endpoint} permits
coefficient extraction term by term. The specializations $k=1,2,3,4$
give \eqref{eq:endpoint-R1}--\eqref{eq:endpoint-R4}; in particular,
the terms containing $\pi^2/6$ come from the cubic term of
$\sin(\pi z)/\pi$, while $c_{j,0}=(-1)^j$ cancels the outer
$(-1)^j$ in the highest moment.
\end{proof}

For example, the first semiprime coefficient is
\begin{equation}\label{eq:semiprime-SD}
 R_2(x)=\frac{x}{\log x}\{\log\log x+M\}
       +O\!\left(\frac{x(1+\log\log x)}{(\log x)^2}\right),\qquad
 M=\gamma_E-\sum_{m\ge2}\frac{P(m)}m.
\end{equation}
Indeed $G(1,z)=1-z\sum_{m\ge2}P(m)/m+O(z^2)$ and
$1/\Gamma(z)=z+\gamma_Ez^2+O(z^3)$. Thus $M$ is the classical
prime Mertens constant defined above, not a new coefficient of this construction.

The expansion does not uniquely specify the function. If $R_k^{>1/2}$
denotes the integral in \eqref{eq:Rk} restricted to $1/2<\sigma<1$,
then Lemma~\ref{lem:cut} gives
\begin{equation}\label{eq:belowhalf}
 R_k(x)-R_k^{>1/2}(x)=O_k(\sqrt{x}).
\end{equation}
At $k=1$,
$R_1^{>1/2}(x)=\li(x)-\li(\sqrt{x})$, whereas $R_1(x)=R(x)$.
These functions and $\li(x)$ have the same expansion to every fixed
inverse-logarithmic order.

\paragraph{What differs from Selberg--Delange.}
The distinction is between a fixed-order asymptotic expansion and a
prescribed approximating function. Proposition~\ref{prop:SD} gives
exactly the classical Selberg--Delange coefficients. Definition
\eqref{eq:Rk} additionally specifies the full positive real cut,
including the lower-cut term in \eqref{eq:belowhalf}, which every
fixed inverse-logarithmic truncation leaves unresolved. Refined
Selberg--Delange contour analyses can also retain secondary terms;
the distinction is not a limitation of that method. The subsequent
results concern the discrepancy of this prescribed function: RH gives
its leading zero spectrum, simplicity aligns consecutive spectra, and
\eqref{eq:moment} supplies the second-order law. These conditional
averaged statements are not new Selberg--Delange coefficients.
The normalization $R_1=R$ alone gives no uniqueness theorem, and the
definition gives neither positivity of $\Delta_k$ nor optimal numerical
accuracy, including against an optimally truncated expansion.

\section{Perron inversion and the arithmetic discrepancy}\label{sec:perron}
The preceding construction concerns a smooth function. Perron's formula
specifies what must be added to recover the count. For $c>1$, set
\[
 C_{k,T}(x)=\frac1{2\pi\ii}\int_{c-\ii T}^{c+\ii T}
                          F_k(s)\frac{x^s}{s}\dd s.
\]
Then $C_{k,T}(x)\to N_k^*(x)$, where
$N_k^*(x)=N_k(x)-\tfrac12\mathbf1_{\{x\in\mathbb N,\,\Omega(x)=k\}}$
is the half-weight count \cite{Titchmarsh}.

For a finite deformation choose $\varepsilon\in(0,1)$ and $T$
so that the boundary of $\varepsilon\le\Re s\le c$, $|\Im s|\le T$
contains no singularity. Only finitely many terms of the M\"obius
series in \eqref{eq:generator} can be singular there: their locations
are $1/(nj)$ and $\rho/(nj)$ with $nj\varepsilon\le1$ and
$1\le j\le k$. Group coincident singularities and slit the rectangle
horizontally from them to its left edge, merging overlapping slits.
Continue the logarithms from $\Re s>1$.

Let $Z_{k,\varepsilon,T}$ be the sum of the lower-minus-upper bank
integrals of $F_k(s)x^s/s$ on the nonreal slits, divided by $2\pi\ii$;
the lower banks run from left to right. Let $E_{k,\varepsilon,T}$
be minus the integral, with the same normalization, along the outer
path from $c+\ii T$ across the top, down the left edge, and across
the bottom to $c-\ii T$, with the inherited branches. Cauchy's theorem gives
\begin{equation}\label{eq:finiteperron}
 C_{k,T}(x)=\int_\varepsilon^1\Delta_k(\sigma)
                       \frac{x^\sigma}{\sigma}\dd\sigma
              +Z_{k,\varepsilon,T}(x)+E_{k,\varepsilon,T}(x).
\end{equation}
The endpoint arcs vanish since their integrals are
$O(\eta(1+|\log\eta|)^k)$ as the radius $\eta\downarrow0$.
First taking $T\to\infty$ through admissible heights, then
$\varepsilon\downarrow0$, gives the unconditional identity
\begin{equation}\label{eq:exactcount}
 N_k^*(x)=R_k(x)+
       \lim_{\varepsilon\downarrow0}\lim_{T\to\infty}
         \{Z_{k,\varepsilon,T}(x)+E_{k,\varepsilon,T}(x)\}.
\end{equation}
The combined limit exists by \eqref{eq:finiteperron}.
Separate convergence of its two terms, or an exact formula involving
zeros alone, is not asserted.

We now study $D_k=R_k-N_k$. The fixed-order remainder estimate in \eqref{eq:SD}
is much larger than $\sqrt{x}(\log\log x)^{k-1}/\log x$; that
expansion alone does not justify replacing the full $R_k$ in the
following normalization. Put
\begin{equation}\label{eq:normalization}
 A_k(y)=(-1)^{k-1}(k-1)!
       \frac{ye^{-y/2}}{(\log y)^{k-1}}D_k(e^y),\qquad y\ge2,
\end{equation}
and
$M_k(Y)=Y^{-1}\int_2^Y|A_{k+1}(y)-A_k(y)|^2\dd y$.
The bounded difference $N_k-N_k^*$ does not affect these mean-square
limits. Under RH the contour can be fixed at $\Re s=2/5$ and its
boundary error becomes negligible after normalization, as shown next.

\section{The first-order spectrum}\label{sec:first}
Under RH, let $\mathcal Z$ be the set of distinct nontrivial zeros
$\rho=\tfrac12+\ii\gamma$, and let $m_\rho$ be the multiplicity.
Every zero sum in this section includes both signs of $\gamma$.
All ranks are fixed. The functions $A_k$ and $M_k$ have the normalization
in~\eqref{eq:normalization}; the complete real-cut function, not a finite
Selberg--Delange truncation, is used in that definition.

\begin{theorem}\label{thm:first}
Assume RH. For every fixed $k\ge1$ and $T\ge2$,
\begin{equation}\label{eq:first}
 \limsup_{Y\to\infty}\frac1Y\int_2^Y
 \left|A_k(y)-\sum_{\substack{\rho\in\mathcal Z\\|\gamma|\le T}}
       \frac{m_\rho^k}{\rho}e^{\ii\gamma y}\right|^2\dd y
 \ll_k\frac{\log^{2k}(T+3)}{T}.
\end{equation}
If the zeros are simple, the right side may be replaced by
$O_k(\log^2(T+3)/T)$. Consequently,
\begin{equation}\label{eq:defect}
 \lim_{Y\to\infty}M_k(Y)
  =\sum_{\rho\in\mathcal Z}
       \frac{m_\rho^{2k}(m_\rho-1)^2}{|\rho|^2}.
\end{equation}
Under RH, vanishing adjacent mean square for any one fixed pair is
therefore equivalent to simplicity of all nontrivial zeros.
\end{theorem}

With simple zeros the adjacent conclusion is equivalently
\begin{equation}\label{eq:raw}
 \lim_{X\to\infty}\frac1{\log X}\int_{e^2}^{X}
 \left|\frac{\log x}{\sqrt{x}(\log\log x)^k}
 \left(D_{k+1}(x)+\frac{\log\log x}{k}D_k(x)\right)\right|^2
 \frac{\dd x}{x}=0.
\end{equation}
It is not a pointwise ratio or a pointwise sign assertion. We prove the
uniform contour and mean-square estimates needed to pass from local
zero contributions to the actual arithmetic error.

\subsection{Uniform contour and local logarithm bounds}
Put $a=2/5$, $d=1/10$ and $l_\rho=\log(|\gamma|+3)$. The usual
local zero count, with multiplicities, gives
\begin{equation}\label{eq:unitcount}
 \sum_{n\le\gamma<n+1}m_\rho\ll\log(|n|+3),\qquad n\in\mathbb Z.
\end{equation}
We use horizontal continuation from $\Re s>1$, with conjugate paths
at conjugate zeros. On the two banks of the cut at $\rho$ write
\begin{equation}\label{eq:local}
 P(\rho-u\pm\ii0)=m_\rho(\log u\pm\ii\pi)+H_\rho(u),
 \qquad 0<u\le d.
\end{equation}
The same $H_\rho(u)$ occurs on both banks. Define
\[
 n_\rho=\sum_{0<|\gamma'-\gamma|\le1}m_{\rho'},\qquad
 L_\rho^*=\sum_{0<|\gamma'-\gamma|\le1}
          m_{\rho'}\log\frac1{|\gamma'-\gamma|}.
\]
Empty sums are zero; coincident zeros have already been grouped in
$m_\rho$.

\begin{lemma}[Local bounds and good heights]\label{lem:localbounds}
Under RH, $n_\rho\ll l_\rho$ and
\begin{equation}\label{eq:Hbound}
 |H_\rho(u)|\ll l_\rho+
        \min\{n_\rho\log(1/u),L_\rho^*\},\qquad 0<u\le d.
\end{equation}
Each $H_\rho$ has a finite limit at $u=0$. Uniformly on these cuts,
$P(2s)\ll l_\rho$, while $P(js)\ll_j1$ for $j\ge3$.
For every sufficiently large $H$ there is a height $T\in[H,H+1]$
such that, simultaneously for every fixed rank $r$,
\[
 \sup_{a\le\sigma\le2}|F_r(\sigma\pm\ii T)|\ll_r\log^r(T+3).
\]
On the left boundary, using either bank at a cut endpoint,
$F_r(a+\ii t)\ll_r\log^r(|t|+3)$.
\end{lemma}
\begin{proof}
The local logarithmic-derivative formula and zero count are standard
\cite[Chapter IX]{Titchmarsh}; their use in almost-prime contours is
also recorded in \cite[Lemmas 9--11]{Meng}. At a large ordinate $t$,
integration of
\[
 \frac{\zeta'}{\zeta}(s)=
 \sum_{|\gamma'-t|\le1}\frac{m_{\rho'}}{s-\rho'}+O(\log(|t|+3))
\]
from $2+\ii t$ along a horizontal path determines the logarithm,
including its branch. Boundary choices at distance exactly $1$ alter
the error by $O(\log(|t|+3))$ only. Under RH,
\[
 P(s)=\log\zeta(s)-\tfrac12\log\zeta(2s)+G_0(s),
 \qquad \Re s\ge a,
\]
where the series defining $G_0$ is absolutely and uniformly bounded.
For a nonreal zero cut, $\log\zeta(2s)=O(l_\rho)$: its argument has
real part at least $4/5$, a fixed distance to the right of all zeros.
The terms with scaling index at least three are absolutely convergent.
After removing the central logarithm, the other local logarithms have
arguments $-u+\ii(\gamma-\gamma')$. Their imaginary parts are bounded,
and their absolute values are at most
\[
 C+\min\{\log(1/u),\log(1/|\gamma-\gamma'|)\}.
\]
Summation proves~\eqref{eq:Hbound}. For a fixed zero all noncentral
local factors are nonzero at $u=0$, so the remaining part is continuous
there. The assertions about $P(js)$ follow from the same formula for
$j=2$ and from absolute convergence for $j\ge3$.

For good heights, integrate over $H\le t\le H+1$ the nonnegative sum
\[
 \sum_{|\gamma'-t|\le1}m_{\rho'}\log^+\frac1{|t-\gamma'|}.
\]
Its integral is $O(\log(H+3))$ by~\eqref{eq:unitcount} and
$\int_{-1}^{1}\log(1/|v|)\dd v<\infty$. Choose a nonsingular $T$
where that bound holds. The integrated local formula gives
$\log\zeta(\sigma+\ii T)=O(\log(T+3))$ throughout $a\le\sigma\le2$,
with the branch just specified. The other prime-zeta factors obey the
same or stronger bounds. The partition polynomial for $F_r$ proves
the horizontal assertion. At $\sigma=a$ all zero factors have distance
at least $d$; the same argument proves the vertical bound. Bounded
ordinates contribute only a fixed constant, with one-sided values at
the real cut.
\end{proof}

Choose $T_U\in[e^{2U},e^{2U}+1]$ as in the lemma for each dyadic
$U=2^j\ge2$, and set $T(y)=T_U$ when $U\le y<2U$.
Thus $e^y\le T(y)\le2e^{2y}$, and the cutoff is constant on each block.
It is the same cutoff for all the finitely many ranks being compared.

\begin{lemma}[Uniform finite contour]\label{lem:uniformcontour}
For $y\in[U,2U)$ and fixed $r$, define
\[
 C_{r,\rho}(y)=\frac1{2\pi\ii}\int_0^d e^{-uy}
 \frac{F_r(\rho-u-\ii0)-F_r(\rho-u+\ii0)}{\rho-u}\dd u.
\]
Under RH,
\begin{equation}\label{eq:perron}
 N_r(e^y)-R_r(e^y)=e^{y/2}
  \sum_{|\gamma|\le T_U}e^{\ii\gamma y}C_{r,\rho}(y)
   +O_r\bigl(e^{ay}(1+y)^{r+1}+1\bigr).
\end{equation}
The error, after the normalization for $A_r$ and multiplication by
any fixed power of $\log y$, tends to zero exponentially.
\end{lemma}
\begin{proof}
The coefficients of $F_r$ lie in $\{0,1\}$. The truncated Perron
kernel therefore gives, uniformly for real $x\ge e^2$ and
$c=1+1/\log x$,
\[
 N_r^*(x)=\frac1{2\pi\ii}\int_{c-\ii T}^{c+\ii T}
       F_r(s)\frac{x^s}{s}\dd s+O(x\log x/T+1).
\]
For completeness, the kernel error is bounded by the sum of
$(x/n)^c\min(1,(T|\log(x/n)|)^{-1})$ away from $n=x$.
The at most two closest integers contribute $O(1)$; the other
integers in $[x/2,2x]$ give $O(x\log x/T)$ by a harmonic sum,
and the two outer ranges give the same bound. At an integer the
half-weight convention handles the central term. Changing to $N_r$
costs at most $1/2$.

In $a\le\Re s\le c$, RH leaves only the nonreal points $\rho$ and
the real points $1$ and $1/2$ as branch points. Remove the disjoint
horizontal nonreal cuts $[a+\ii\gamma,\rho]$ and the single combined
real cut $[a,1]$. Multiple zeros are treated as one branch point.
At each fixed $T$ there are finitely many cuts. Take all endpoint
radii to zero before allowing $T$ to grow. Their integrals vanish as
$O_\rho(\eta(1+|\log\eta|)^r)$; no radius uniform in zero gaps is
needed. The cut banks are the horizontal continuations of
Lemma~\ref{lem:localbounds}, so their values do not depend on $T$.

The horizontal sides cost $O_r(x\log^r(T+3)/T)$ and the left side
costs $O_r(x^a\log^{r+1}(T+3))$. The real banks give exactly
$\int_a^1\Delta_r(\sigma)x^\sigma\dd\sigma/\sigma$.
Lemma~\ref{lem:cut} bounds the omitted real-cut portion by $O_r(x^a)$.
The nonreal banks give the displayed $C_{r,\rho}$. Substitution of
$T=T_U$ proves~\eqref{eq:perron}. In particular the boundary estimate
is independent of the number of small endpoint circles.
\end{proof}

\subsection{A Gaussian--Laplace mean-square lemma}
For a finite set of ordinates and complex coefficients define
\[
 \mathcal B(c)^2=\sum_{n\in\mathbb Z}
      \left(\sum_{n\le\gamma<n+1}|c_\rho|\right)^2,
 \qquad
 \|f\|_U^2=\frac1U\int_U^{2U}|f(y)|^2\dd y.
\]
The same notation is used for an infinite coefficient sequence when
the displayed sum converges.

\begin{lemma}[Uniform block and Laplace bounds]\label{lem:gaussianlaplace}
For $U\ge1$ and every finite frequency set,
\begin{equation}\label{eq:gaussian}
 \left\|\sum_\rho c_\rho e^{\ii\gamma y}\right\|_U
       \le C\mathcal B(c).
\end{equation}
For $U\ge2$, fixed real $p$, and measurable coefficient functions
$a_\rho(u)$ on $(0,d)$ for which the right side is finite,
\begin{align}\label{eq:laplaceblock}
 &\left\|\frac{y}{(\log y)^p}\sum_\rho e^{\ii\gamma y}
       \int_0^d e^{-uy}a_\rho(u)\dd u\right\|_U\notag\\
 &\hspace{10mm}\le C_p\int_0^{dU}e^{-v}
       \mathcal B\left(\left\{
            \frac{a_\rho(v/U)}{(\log U)^p}\right\}_\rho\right)\dd v.
\end{align}
Both constants are independent of the number, spacing and maximum
height of the frequencies. On zeta zeros,
\begin{equation}\label{eq:binweighted}
 \mathcal B(c)^2\ll\sum_\rho l_\rho|c_\rho|^2.
\end{equation}
\end{lemma}
\begin{proof}
Majorize $1_{[U,2U]}(y)$ by
$e^{1/4}\exp(-((y-3U/2)/U)^2)$ and expand the square. The normalized
Fourier kernel has absolute value at most
$C\exp(-U^2(\gamma-\gamma')^2/4)$. If the frequencies lie in bins
$n,m$, this is bounded, since $U\ge1$, by
$C\exp(-((|n-m|-1)_+)^2/4)$. This sequence in $n-m$ is summable.
The $\ell^2$ convolution inequality proves~\eqref{eq:gaussian}.
No minimum spacing has been used, even within one bin.

Minkowski's integral inequality, $y\le2U$,
$e^{-uy}\le e^{-Uu}$ and comparability of $\log y$ and $\log U$ give
an upper bound
\[
 \frac{C_pU}{(\log U)^p}\int_0^d e^{-Uu}
       \mathcal B(\{a_\rho(u)\})\dd u.
\]
Set $v=Uu$ to obtain~\eqref{eq:laplaceblock}.
Finally apply Cauchy--Schwarz in each unit bin and use
\eqref{eq:unitcount}. The quantities $l_\rho$ are comparable within
one bin, proving~\eqref{eq:binweighted}.
\end{proof}

We will repeatedly use the following precise dominated-convergence
consequence. Suppose $c_{\rho,U}(v)\to0$ for every fixed $\rho$ and
$v>0$, and
\begin{equation}\label{eq:domcriterion}
 |c_{\rho,U}(v)|\le d_\rho(1+|\log v|)^M,
 \qquad \sum_\rho l_\rho d_\rho^2<\infty.
\end{equation}
Set the coefficient to zero when $v>dU$. By~\eqref{eq:binweighted}
and two applications of dominated convergence,
\[
 \int_0^{dU}e^{-v}\mathcal B(\{c_{\rho,U}(v)\})\dd v\longrightarrow0.
\]
The conclusion holds uniformly when any finite subset of the zeros
is retained, including $|\gamma|\le T_U$. For a height tail the same
bound uses $\sum_{|\gamma|>T}l_\rho d_\rho^2$ instead. This is how
height uniformity is obtained; coefficients are not frozen before
estimating the high-zero terms.

\subsection{Completion of the first-order proof}
\begin{proof}[Proof of Theorem~\ref{thm:first}]
The partition formula writes $F_r$ as its top term $P^r/r!$ plus
finitely many terms
\begin{equation}\label{eq:partitionaudit}
 c_\nu P(s)^{r_\nu}\prod_{j\ge2}P(js)^{\nu_j},
 \quad r_\nu+\sum_{j\ge2}j\nu_j=r,\quad r_\nu\le r-2.
\end{equation}
Put $q_\nu=\sum_{j\ge2}\nu_j\le(r-r_\nu)/2$. On a cut the product
of scaled factors is analytic and bounded by $C_r l_\rho^{q_\nu}$.
Expanding~\eqref{eq:local}, every nonzero jump is a constant multiple of
\begin{equation}\label{eq:jumpterm}
 \frac{m_\rho^bJ_b(u)H_\rho(u)^c g_\rho(u)}{\rho-u},
 \quad b\ge1,\quad b+c=r_\nu,\quad |g_\rho(u)|\le C_rl_\rho^q,
\end{equation}
where the top term has $r_\nu=r$, $q=0$, and
\[
 J_b(u)=\frac{(\log u-\ii\pi)^b-(\log u+\ii\pi)^b}{2\pi\ii},
 \qquad |J_b(u)|\le C_b(1+|\log u|)^{b-1}.
\]
A term with $b=0$ has no jump. The coefficients $g_\rho(u)$ are
independent of $y$ and of the height cutoff.

Consider first every term except $b=r,c=q=0$. Apply
\eqref{eq:laplaceblock} with $p=r-1$. Write $\ell_U=\log U$ and
substitute $u=v/U$. Since $m_\rho\ll l_\rho$ and
$|H_\rho(u)|\ll l_\rho(1+\log(1/u))$, its normalized coefficient
is bounded by
\begin{equation}\label{eq:firstenvelope}
 \frac{C_r l_\rho^{b+c+q}}{|\rho|}
       (1+|\log v|)^{b+c-1},\qquad 0<v<dU.
\end{equation}
Here $b+c\le r$; all remaining powers of $\ell_U$ are bounded.
For fixed $\rho,v$, the regular factors tend to finite values, and
$b<r$ makes the same coefficient tend to zero, with a factor
$\ell_U^{b-r}$. For every fixed $A$,
$\sum_\rho l_\rho^A/|\rho|^2<\infty$ by~\eqref{eq:unitcount}.
Thus~\eqref{eq:domcriterion} applies to~\eqref{eq:firstenvelope}.
All mixed and lower-partition terms are $o(1)$ in block mean square,
uniformly for the moving cutoff. This includes arbitrarily close
but distinct zeros.

For the pure top term replace $(\rho-u)^{-1}$ by $\rho^{-1}$.
The difference is $O(u/|\rho|^2)$. After $u=v/U$, the block bound is
$O_r(U^{-1})$, since the corresponding weighted coefficient square
sum converges. Extending its integral from $d$ to infinity costs
$O_r(e^{-dU/2})$ in block norm. Both statements remain true after
multiplication by $\log y$, with errors tending to zero. The exact
Laplace identity is
\begin{equation}\label{eq:purelaplace}
 y\int_0^\infty e^{-uy}J_r(u)\dd u
   =\left.\frac{d^r}{dz^r}
       \frac{e^{-z\log y}}{\Gamma(-z)}\right|_{z=0}.
\end{equation}
Indeed $\sum_{r\ge0}J_r(u)z^r/r!=-u^z\sin(\pi z)/\pi$ and
$-\Gamma(1+z)\sin(\pi z)/\pi=1/\Gamma(-z)$; differentiation near
$z=0$ is dominated. These are the usual Hankel--Gamma identities
\cite{DLMF}. With
\[
 \Pi_r(\ell)=\frac{(-1)^r}{r\ell^{r-1}}
       \left.\frac{d^r}{dz^r}
          \frac{e^{-\ell z}}{\Gamma(-z)}\right|_{z=0},
\]
one has
\begin{equation}\label{eq:Pi}
 \Pi_r(\ell)=1+\frac{(r-1)\gamma_E}{\ell}+O_r(\ell^{-2}).
\end{equation}
It is a finite polynomial in $1/\ell$, with $\Pi_1=1$.
The sign $D_r=R_r-N_r$ and its factorial normalization consequently give
\[
 A_r(y)=\Pi_r(\log y)S_{r,U}(y)+E_r(y),\qquad
 S_{r,U}(y)=\sum_{|\gamma|\le T_U}\frac{m_\rho^r}{\rho}e^{\ii\gamma y},
 \quad \|E_r\|_U=o(1).
\]
The contour error is included here by Lemma~\ref{lem:uniformcontour}.

For an arbitrary finite subset of $|\gamma|>T$,
\begin{equation}\label{eq:tail}
 \mathcal B\left(\left\{\frac{m_\rho^r}{\rho}
                     1_{\{|\gamma|>T\}}\right\}\right)^2
 \ll_r\frac{\log^{2r}(T+3)}T.
\end{equation}
To see this, in one bin use
$\sum m_\rho^r\le(\sum m_\rho)^r\ll_r\log^r(|n|+3)$ and then
sum $\log^{2r}(|n|+3)/(1+|n|)^2$ over the tail. With simple zeros
the bin sum is $O(\log(|n|+3))$, independently of $r$.
In particular $\|S_{r,U}\|_U$ is bounded uniformly, so
$\Pi_r-1=o(1)$ causes no difficulty. For fixed $T$, eventually
$T_U>T$; apply~\eqref{eq:gaussian} and~\eqref{eq:tail} to the
moving part above $T$.

If block norms tend to zero, their full Ces\`aro mean square tends
to zero: split $[2,Y]$ into dyadic blocks, bound a final partial block
by its complete block, and divide by $Y$. The same argument transfers
a uniform squared block bound to a uniform limsup bound. Finite initial
blocks disappear. This proves~\eqref{eq:first} without interchanging
the infinite zero sum and the limit in $Y$.

Distinct-frequency orthogonality for fixed truncations, followed by
\eqref{eq:first}, gives
\[
 \lim_{Y\to\infty}\frac1Y\int_2^Y A_j(y)A_k(y)\dd y
   =\sum_{\rho\in\mathcal Z}\frac{m_\rho^{j+k}}{|\rho|^2}.
\]
Apply this to the difference to obtain~\eqref{eq:defect}.
Each summand in that defect is nonnegative and vanishes precisely
when $m_\rho=1$. No linear-independence assumption on the ordinates
is needed.
\end{proof}

\paragraph{Leading covariance.}
Under RH and simplicity the leading covariance for every fixed pair
of ranks is
\[
 V=\sum_\rho\frac1{|\rho|^2}=2+\gamma_E-\log(4\pi)>0.
\]
The evaluation follows from the Hadamard product for $\xi$ at $0$ and
$1$ \cite{Titchmarsh}. Thus a fixed finite family has a rank-one leading
covariance matrix. Removing $(-1)^{k-1}$ from the normalization
restores the alternating signs.

\section{The second-order spectrum}\label{sec:second}
Assume RH and simplicity. Now $L_\rho=L_\rho^*$, and for a fixed pair
$(k,k+1)$ impose
\begin{equation}\tag{$\mathrm H_k$}\label{eq:moment}
 \sum_\rho\frac{l_\rho^{2k-1}(l_\rho+L_\rho)^2}{|\rho|^2}<\infty.
\end{equation}
Simplicity makes each term finite; it does not imply convergence of
the series. Define the branch-consistent finite part
\begin{equation}\label{eq:coefficients}
 h_\rho=\lim_{v\downarrow0}\{P(\rho+v)-\log v\}=H_\rho(0),
 \qquad b_\rho=\frac{\gamma_E-h_\rho}{\rho},
 \qquad W=\sum_\rho|b_\rho|^2.
\end{equation}
Independent choices of principal logarithms of $\zeta'(\rho)$ are
not used.

\begin{theorem}\label{thm:second}
Under RH, simplicity and~\eqref{eq:moment}, $0<W<\infty$ and
\begin{equation}\label{eq:second}
 \lim_{T\to\infty}\limsup_{Y\to\infty}\frac1Y\int_2^Y
 \left|\log y\,(A_{k+1}(y)-A_k(y))
       -\sum_{|\gamma|\le T}b_\rho e^{\ii\gamma y}\right|^2\dd y=0.
\end{equation}
Consequently,
\begin{equation}\label{eq:rate}
 M_k(Y)\sim\frac{W}{(\log Y)^2}.
\end{equation}
The constant is the same for every fixed pair for which its required
moment converges.
\end{theorem}

\subsection{Second-order domination with the exact moment}
Put
\[
 S_U(y)=\sum_{|\gamma|\le T_U}\frac{e^{\ii\gamma y}}\rho,
 \qquad B_U(y)=\sum_{|\gamma|\le T_U}b_\rho e^{\ii\gamma y}.
\]
\begin{lemma}[Second-order block expansion]\label{lem:secondblock}
Under the hypotheses of Theorem~\ref{thm:second}, for every
$1\le r\le k+1$,
\begin{equation}\label{eq:secondblock}
 \left\|\log y\{A_r(y)-S_U(y)\}-(r-1)B_U(y)\right\|_U\longrightarrow0
 \quad(U=2^j\longrightarrow\infty).
\end{equation}
\end{lemma}
\begin{proof}
Use the same partition and jump terms as in~\eqref{eq:jumpterm}, now
with $m_\rho=1$. Set
$A_\rho=1+L_\rho/l_\rho$ and
$\tau_U(v)=1+\log U+|\log v|$. Lemma~\ref{lem:localbounds} gives
\begin{equation}\label{eq:secondH}
 |H_\rho(v/U)|\le C l_\rho\min\{\tau_U(v),A_\rho\},
 \qquad 0<v<dU.
\end{equation}
For $c\ge1$, the elementary inequality
$\min(\tau,A)^c\le A\tau^{c-1}$ implies that, after the second-order
normalization $p=r-2$ in~\eqref{eq:laplaceblock}, a jump term is
bounded by
\begin{equation}\label{eq:secondenvelope}
 \frac{C_r l_\rho^{q+c-1}(l_\rho+L_\rho)}{|\rho|}
 \frac{\tau_U(v)^{b+c-2}}{(\log U)^{r-2}}.
\end{equation}
For a top-partition term $q=0$ and $b+c=r$; for a lower term
$b+c=r_\nu\le r-2$ and $q\le(r-r_\nu)/2$. Since $b\ge1$,
in either case
\[
 q+c-1\le r-2\le k-1.
\]
(The lower-partition bound is in fact at most $r-3$.) Thus
\eqref{eq:secondenvelope} is bounded by
\begin{equation}\label{eq:exactHdom}
 \frac{C_k l_\rho^{k-1}(l_\rho+L_\rho)}{|\rho|}
              (1+|\log v|)^{k-1}.
\end{equation}
The unit-bin square sum of this envelope is controlled by exactly
\eqref{eq:moment}, because its square acquires the one additional
factor $l_\rho$ in~\eqref{eq:binweighted}.

For the top terms with $c\ge2$, a fixed zero and fixed $v$ give
$H_\rho(v/U)=O_\rho(1)$, so the normalized coefficient has the
factor $(\log U)^{1-c}\to0$. Every lower-partition term with $c\ge1$
also tends to zero, since $b\le r-2$. Dominated convergence with
\eqref{eq:exactHdom} eliminates these terms in block norm, uniformly
in the height cutoff. Lower-partition terms with $c=0$ have
$b\le r-2$, a factor $(\log U)^{b-r+1}\to0$, and an envelope
$C_r l_\rho^q(1+|\log v|)^{b-1}/|\rho|$ with a summable weighted
square. No gap moment is needed for this last case.

Two contributions remain. The pure top term was calculated in
\eqref{eq:purelaplace}--\eqref{eq:Pi}; its denominator and endpoint
errors still vanish after multiplication by $\log y$. Since
$\|S_U\|_U$ is uniformly bounded, it contributes
\[
 \log y\,S_U(y)+(r-1)\gamma_E S_U(y)+o_{L^2(U,2U)}(1)
\]
to $\log y\,A_r(y)$, where $o_{L^2(U,2U)}(1)$ means that the
normalized block norm tends to zero.

For $r\ge2$, the top jump with $c=1,b=r-1$ is
$J_{r-1}(u)H_\rho(u)/(r-1)!$ before division by $\rho-u$. Inside its Laplace integral
replace $H_\rho(u)/(\rho-u)$ by $h_\rho/\rho$.
This replacement is valid globally in block norm, not just at finitely
many zeros: the difference tends to zero for each fixed $\rho,v$
and is dominated by
\[
 \frac{C_r(l_\rho+L_\rho)}{|\rho|}
                  (1+|\log v|)^{r-2}.
\]
Its weighted square is summable under~\eqref{eq:moment}.
Apply~\eqref{eq:domcriterion}, retaining the original coefficient
functions until this domination has been established.
Formula~\eqref{eq:purelaplace} for $J_{r-1}$ then shows that this term
contributes
\[
 -(r-1)\Pi_{r-1}(\log y)
          \sum_{|\gamma|\le T_U}\frac{h_\rho}{\rho}e^{\ii\gamma y}
 =-(r-1)\sum_{|\gamma|\le T_U}\frac{h_\rho}{\rho}e^{\ii\gamma y}
       +o_{L^2(U,2U)}(1).
\]
The last equality uses~\eqref{eq:gaussian},
$|h_\rho|\ll l_\rho+L_\rho$ and the same moment. When $r=1$ this
term is absent. Combining the two surviving terms proves
\eqref{eq:secondblock}. All other partition terms and the contour
error have been accounted for.
\end{proof}

\subsection{The limiting spectrum, the rate and positivity}
\begin{lemma}[Logarithmic weight transfer]\label{lem:weighttransfer}
If $B$ is locally square integrable on $[2,\infty)$ and
$\int_2^Y|B(y)|^2\dd y=WY+o(Y)$, then
\[
 \frac1Y\int_2^Y\frac{|B(y)|^2}{(\log y)^2}\dd y
       =\frac{W}{(\log Y)^2}+o((\log Y)^{-2}).
\]
\end{lemma}
\begin{proof}
Write $F(t)=\int_2^t|B(y)|^2\dd y=Wt+r(t)$, $r(t)=o(t)$.
Stieltjes partial summation gives
\[
 \int_2^Y\frac{\dd F(t)}{(\log t)^2}
   =\frac{F(Y)}{(\log Y)^2}
       +2\int_2^Y\frac{F(t)}{t(\log t)^3}\dd t.
\]
The second term from $Wt$ is $O(Y/(\log Y)^3)+O(1)$. For the
$r(t)$ terms, split at a fixed point beyond which $|r(t)|\le\varepsilon t$;
the fixed initial part is negligible, and then let $\varepsilon\downarrow0$.
\end{proof}

\begin{proof}[Proof of Theorem~\ref{thm:second}]
Subtract~\eqref{eq:secondblock} at $r=k+1$ and $r=k$.
The moving leading spectrum cancels exactly, leaving
\[
 \|\log y\,(A_{k+1}-A_k)-B_U\|_U\longrightarrow0.
\]
For any finite subset of the high zeros,
\begin{equation}\label{eq:btail}
 \mathcal B(\{b_\rho 1_{\{|\gamma|>T\}}\})^2
  \ll\sum_{|\gamma|>T}\frac{l_\rho(l_\rho+L_\rho)^2}{|\rho|^2}
       \longrightarrow0.
\end{equation}
This follows from~\eqref{eq:binweighted} and is dominated by the tail
of~\eqref{eq:moment}. Apply~\eqref{eq:gaussian} block by block, then
pass to the full average exactly as in the first-order proof. First
$Y$ tends to infinity with $T$ fixed; only then does $T$ tend to
infinity. This proves~\eqref{eq:second} and $W<\infty$.
Finite-frequency orthogonality and the tail control give
\[
 \frac1Y\int_2^Y|\log y\,(A_{k+1}-A_k)|^2\dd y\longrightarrow W.
\]
Lemma~\ref{lem:weighttransfer} proves~\eqref{eq:rate} once $W>0$.

More generally, the finite parts $h_\rho$ cannot all equal one real
constant $c$. At a simple zero the real part of the M\"obius identity
is, uniformly in height,
\begin{equation}\label{eq:hderivative}
 \Re h_\rho=\log|\zeta'(\rho)|-\tfrac12\log|\zeta(2\rho)|+O(1).
\end{equation}
Indeed the terms with scaling index at least three converge absolutely
at $\Re\rho=1/2$ and have a uniformly bounded sum. If $h_\rho=c$,
then $|\zeta'(\rho)|^2\ll_c|\zeta(1+2\ii\gamma)|\ll_c l_\rho$,
using the standard bound on the line $\Re s=1$.
The sum of these squares up to height $T$ would be $O_c(T\log^2T)$.
But Gonek's RH theorem \cite[Corollary 1, $\mu=1$, $\alpha=0$]{Gonek}
gives
\[
 \sum_{0<\gamma\le T}|\zeta'(\rho)|^2
   =\frac{T\log^4T}{24\pi}+O(T\log^3T).
\]
This is a contradiction. In particular $W=0$ would require
$h_\rho=\gamma_E$ at every zero, so $W>0$.
\end{proof}

\paragraph{Sufficient hypotheses for the extra moment.}
A polynomial lower bound on neighboring gaps implies
\eqref{eq:moment}. So does a bound
$|\gamma-\gamma'|\ge\exp(-C(|\gamma|+3)^\alpha)$ for
$0<|\gamma-\gamma'|\le1$ with fixed $\alpha<1/2$:
then $L_\rho\ll l_\rho(|\gamma|+3)^\alpha$, and the resulting
series converges by unit-bin counting. Neither gap bound is being
deduced from RH and simplicity.

The Weak Mertens Conjecture is another sufficient hypothesis. Its
standard consequences include RH, simplicity and
$\sum_\rho|\rho\zeta'(\rho)|^{-2}<\infty$
\cite[Section 1.2]{BFM}. Taking real parts in the integrated local
formula and using~\eqref{eq:hderivative} gives
$L_\rho\ll l_\rho+\log^+(1/|\zeta'(\rho)|)$.
For every fixed $A$, the inequality $(\log^+t)^2\ll t$ and
Cauchy--Schwarz give
\[
 \sum_\rho\frac{l_\rho^A(\log^+(1/|\zeta'(\rho)|))^2}{|\rho|^2}
 \ll\left(\sum_\rho\frac1{|\rho\zeta'(\rho)|^2}\right)^{1/2}
    \left(\sum_\rho\frac{l_\rho^{2A}}{|\rho|^2}\right)^{1/2}<\infty.
\]
Thus it supplies the required moment for every fixed rank. This is a
sufficient additional hypothesis, not a proof of the moment from RH.

\paragraph{Acknowledgements.}
The author is grateful to the MathOverflow contributor Lucia for the
exposition of the real-axis slit-contour approach in \cite{LuciaMO}.
OpenAI ChatGPT assisted with derivations, checks and manuscript preparation.
The author is responsible for the arguments and references.

\begingroup\small
\begin{thebibliography}{99}
\setlength{\itemsep}{2pt}
\bibitem{SD} R. de la Bret\`eche and G. Tenenbaum,
\emph{Remarks on the Selberg--Delange method},
Acta Arith. \textbf{200} (2021), 349--369.
\bibitem{CrisanErban} D. Cri\c{s}an and R. Erban,
\emph{On the counting function of semiprimes},
Integers \textbf{21} (2021), A122.
\bibitem{Titchmarsh} E. C. Titchmarsh,
\emph{The Theory of the Riemann Zeta-Function}, 2nd ed.,
revised by D. R. Heath-Brown, Oxford University Press, 1986.
\bibitem{BFM} H. M. Bui, A. Florea and M. B. Milinovich,
\emph{Negative discrete moments of the derivative of the Riemann zeta-function},
Bull. London Math. Soc. \textbf{56} (2024), 2680--2703.
\href{https://doi.org/10.1112/blms.13092}{doi:10.1112/blms.13092}.
\bibitem{Froberg} C.-E. Fr\"oberg, \emph{On the prime zeta function},
BIT \textbf{8} (1968), 187--202.
\bibitem{Gonek} S. M. Gonek, \emph{Mean values of the Riemann zeta-function
and its derivatives}, Invent. Math. \textbf{75} (1984), 123--141, Corollary~1.
\bibitem{Mathar} R. J. Mathar, \emph{Series of reciprocal powers of
$k$-almost primes}, arXiv:0803.0900v3 (2009).
\bibitem{LuciaMO} Lucia, answer to ``Asymptotic density of
$k$-almost primes'', MathOverflow, 31 October 2013.
\url{https://mathoverflow.net/questions/35927/asymptotic-density-of-k-almost-primes}.
\bibitem{Meng} X. Meng, \emph{Chebyshev's bias for products of $k$ primes},
Algebra \& Number Theory \textbf{12} (2018), 305--341.
\href{https://doi.org/10.2140/ant.2018.12.305}{doi:10.2140/ant.2018.12.305}.
\bibitem{MV} H. L. Montgomery and R. C. Vaughan,
\emph{Multiplicative Number Theory I: Classical Theory},
Cambridge University Press, 2007.
\bibitem{DLMF} NIST Digital Library of Mathematical Functions,
\S5.7, equation 5.7.3; \S5.9, equation 5.9.2;
\S8.2, equation 8.2.1; and \S25.2, equation 25.2.4.
\url{https://dlmf.nist.gov/5.7}.
\bibitem{BeltonParity} M.~Belton, \emph{Almost-prime error correlations and parity}, companion research manuscript, 6 November 2026.
\end{thebibliography}
\endgroup
\end{document}